\documentclass{article}
\usepackage[a4paper, total={7in, 10in}]{geometry}
\usepackage{amsmath,amsthm,amssymb,mathrsfs,amstext, enumitem,comment,url,hyperref}
\hypersetup{colorlinks,linkcolor={magenta},citecolor={blue}}
\newtheorem{thm}{Theorem}[section]
\newtheorem{lem}[thm]{Lemma}

\usepackage{silence}
\theoremstyle{definition}

\newtheorem{note}[thm]{Note}

\theoremstyle{remark}

\newcommand{\POL}{\mathbf{POL}}

\newcommand{\N}{\mathbb{N}}

\newcommand{\FS}{\operatorname{FS}}
\newcommand*{\rom}[1]{\expandafter\romannumeral #1}

\numberwithin{equation}{section}
\begin{document}

\title{\textbf{An Explicit Ramsey Bound for de Polignac Numbers}}
\author{Sayan Goswami} 
\date{}
\maketitle

\begin{abstract}
A number \( m \) is called a \textbf{de Polignac number}(\( \mathbf{POL} \) in short) if it can be expressed as the difference of infinitely many pairs of consecutive prime numbers. For any given $r\in \N,$ a set $A$ is said to be $\Delta_r^\star$ if for all sets $S$ with $|S|=r$ such that $A\cap \{s-t:s>t\in S\}\neq \emptyset.$ In this article we prove that the set $\mathbf{POL}$ is $\Delta_r^\star$ with the specific computable value of $r,$ where $r=exp(\mathcal{O}(50)).$ Then we prove that there exists a set $E$ with $|E|\leq 2^{50\cdot \prod_{i=1}^{49}p_i} $ (where $(p_i)_i$ is the enumeration of primes) such that $\N=\bigcup_{t\in E} t^{-1}\POL.$

\end{abstract}

\noindent \textbf{Keywords:}  Algebra of the Stone-\v{C}ech compactification of discrete Semigroups, Ultrafilters, Difference set of Primes, Twin prime conjecture, de Polignac numbers, Ramsey theory, $IP$-set,$IP_r$-set, piecewise syndetic set\\
\textbf{Mathematics subject classification:} Primary 37A44, 05D10; Secondary 11E25, 11T30.

\section{Introduction}

In \cite{pi}, Pintz established that the set of differences of prime numbers contains bounded gaps. This result was significantly improved by W.~Huang and X.~Wu \cite{key-3}, who, building upon the Maynard–Tao theorem \cite{key-27}, showed that the difference set of primes is even richer—demonstrating that it possesses bounded gaps not only in the additive sense but also in the multiplicative sense. Further developments in this direction have appeared in \cite{sg1, sg2}.

For any two subsets \( A, B \subseteq \mathbb{N} \), define their difference set by
\[
A - B := \{ a - b : a \in A,\ b \in B,\ a > b \}.
\]
Let \( \mathbb{P} \) denote the set of prime numbers. Then \( \mathbb{P} - \mathbb{P} \) consists of all (positive) differences between pairs of primes. Let $p_n$ denote the $n$-th prime. A number $d \in \N$ is a \emph{de Polignac number} if $d = p_{n+1}-p_n$ for infinitely many $n$, where $p_{n+1},p_n$ are consecutive primes. Denote
\[
\POL = \{ d \in \N : d \text{ occurs infinitely often as a gap between consecutive primes} \}.
\]

For any given $r\in \N,$ a set $A\subseteq \N$ is said to be an $IP_r$ set if there exists a sequence $\langle x_n\rangle_{n=1}^r$ such that $$A=\FS(\langle x_n\rangle_{n=1}^r)=\left\lbrace\sum_{t\in H}x_t:(\neq \emptyset )H\subseteq \{1,2,\ldots ,r\}\right\rbrace.$$ A set $A\subseteq \N$ is said to be an $IP_r^\star$ if $A$ intersects any $IP_r$ set. A set is $A$ is said to be an IP set if there exists an infinite sequence $\langle x_n\rangle_n$ such that $A$ is the set of all finite sums of the sequence $\langle x_n\rangle_n$ with no repetition in sum.

Similarly $A\subseteq \N$ is said to be a $\Delta_r$ set if there exists an finite set $S\subseteq \N$ with $|S|=r$ such that $S-S\subseteq A.$ A set $A\subseteq \N$ is said to be an $\Delta_r^\star$ if $A$ intersects any $\Delta_r$ set.

For any $n\in \N$ and $A\subseteq \N,$ we define $n^{-1}A=\{m\in \N: mn\in A\}.$ A set $A\subseteq \N$ is said to be multiplicative syndetic if there exists a finite set $F\subset \N$ such that $\N=\cup_{t\in F}t^{-1}A.$
It is well known that every $\Delta_r^\star$ set is a multiplicatively syndetic set. In \cite{key-3}, Huang and Wu proved that the set \( \mathbb{P} - \mathbb{P} \) is $\Delta_r^\star$ for some computable number $r.$

The first goal of this article is to prove the following theorem.

\begin{thm}\label{A}
There exists $r\in\N$ such that $\POL$ is a $\Delta_r^*$-set, where $r=exp(\mathcal{O}(50 \log 50)).$
\end{thm}

This doesn't only say that $\POL$ is a $\Delta_r^*$ set but also quantitatively improves it. Now $\POL$ is immediately multiplicatively syndetic. Now our next result improves the size of the translating set.

\begin{thm}\label{B}
    The set $\POL$ is multiplicatively syndetic. In fact there exists a set $E$ with $E\subseteq [1, 2^{\,50\prod_{i=1}^{49}p_i}]$ such that $\N=\bigcup_{t\in E} t^{-1}\POL.$
\end{thm}

\section{Proof of our results}

A finite set \( \mathcal{H} \) of distinct nonnegative integers is called \textbf{admissible} if, for every prime \( p \), it avoids at least one residue class modulo \( p \). Building on the work of Tao and Ziegler \cite{T}, we define an admissible tuple \( \mathcal{H} = \{h_1, h_2, \dots, h_k\} \) as \textbf{prime-generating} if there exist infinitely many natural numbers \( n \) such that all elements of the shifted set \( \{n + h_1, n + h_2, \dots, n + h_k\} \) are simultaneously prime. 

A long-standing conjecture, attributed to Dickson, Hardy, and Littlewood, asserts that every admissible tuple \( \mathcal{H} \) is prime-generating. This remains one of the most challenging open problems in analytic number theory. The following theorem, due to Maynard \cite{key-27} and Tao.

\begin{thm}\cite{key-27}\label{mt}
    For every integer \( m > 2 \), there exists a threshold \( k_m \) such that the following holds: if \( (h_1, h_2, \dots, h_k) \) is an admissible tuple with \( k > k_m \), then the set \( \{n + h_1, \dots, n + h_k\} \) contains at least \( m \) primes for infinitely many values of \( n \in \mathbb{N} \).
\end{thm}

We will use the following theorem of Banks–Freiberg–Turnage-Butterbaugh in our proof.

\begin{thm}\cite{b1}\label{ess}  
    Fix an integer \( m > 2 \), and let \( k_m \) be the threshold given in Theorem \ref{mt}. If \( (h_1, h_2, \dots, h_k) \) is an admissible tuple with \( k > k_m \), then there exists a subset \( \{h'_1, \dots, h'_m\} \subseteq \{h_1, \dots, h_k\} \) such that the set

\[
    \{n + h'_1, \dots, n + h'_m\}
    \]
 consists of \( m \) consecutive primes for infinitely many values of \( n \in \mathbb{N} \).
\end{thm}
\begin{note}\label{n}
    For $m=2$, Polymath gives $k_2\le 50$ \cite[Theorem 3.2 (1)]{key-26}. 
\end{note}

We need the following lemma to prove Theorem \ref{A}. 

\begin{lem}\label{lem:admissible-subset}
For every $k\in\mathbb N$ there exists $R(k)\in\mathbb N$ such that
every finite set $S\subseteq\mathbb N$ with $|S|\geq R(k)$ contains an
admissible subset of cardinality $k$.

More precisely, one may take
\[
R(k)
=
\left\lceil
k\prod_{p\leq k}\frac{p}{p-1}
\right\rceil,
\]
where the product is over all primes $p\leq k$. In particular, from Mertens' theorem $R(k)=\exp(O(k \log k))$ and computable.
\end{lem}

\begin{proof}
Let
\(
S\subseteq\mathbb N,\text{ with } |S|\geq R(k)
\)
 and 
\(
p_1<p_2<\cdots<p_t
\)
be the primes not exceeding $k$. We construct a decreasing sequence
\[
S=S_0\supseteq S_1\supseteq\cdots\supseteq S_t
\]
as follows.

Suppose that $S_{i-1}$ has already been constructed. If the residues
represented by $S_{i-1}$ modulo $p_i$ do not contain all $p_i$
residue classes, put \(S_i=S_{i-1}.\)
Otherwise, by the pigeonhole principle, one residue class modulo $p_i$
contains at most $|S_{i-1}|/p_i$ elements. Choose such a residue class
$r_i\pmod{p_i}$ and put \(S_i= S_{i-1}\setminus(r_i+p_i\mathbb Z).\)

In either case, \(|S_i| \geq |S_{i-1}|\left(1-\frac1{p_i}\right).\) Consequently,
\[
|S_t|
\geq
|S|\prod_{i=1}^t
\left(1-\frac1{p_i}\right)
=
|S|\prod_{p\leq k}
\frac{p-1}{p}.
\]
Since \(|S|
\geq
k\prod_{p\leq k}\frac{p}{p-1},\) we obtain \(|S_t|\geq k.\)

By construction, for every prime $p\leq k$, the set $S_t$ omits
at least one residue class modulo $p$. Let $H\subseteq S_t$ be any
subset with $|H|=k$.

Then for every prime $p\leq k$, we have \(|H\bmod p|\leq p-1.\)
On the other hand, if $p>k$, then \(|H\bmod p|\leq |H|=k<p.\)
Thus $H$ occupies fewer than $p$ residue classes modulo every prime
$p$. Hence $H$ is admissible.
\end{proof}

\begin{proof}[Proof of Theorem \ref{B}]
Let \(k=k_2=50\)
and put \[r=R(50)=
\left\lceil
50\prod_{p\leq 50}\frac{p}{p-1}
\right\rceil=\exp(O(50 \log 50)).\]
Let $S\subseteq\mathbb N$ be arbitrary with \(|S|=r.\)
By Lemma~\ref{lem:admissible-subset}, there is an admissible set \(H=\{b_1,\ldots,b_k\}\subseteq S.\) Applying Theorem~\ref{ess} with $m=2$, we obtain two distinct elements $b_i,b_j\in H$ such that, for infinitely many $n\in\mathbb N$, \(n+b_i, n+b_j\) are consecutive primes.

Without loss of generality, suppose that $b_j>b_i$. Then $b_j-b_i$ occurs infinitely often as a difference of consecutive primes. Therefore $b_j-b_i\in{\rm POL}.$ 
Since $b_i,b_j\in S$, we also have $b_j-b_i\in\Delta(S).$

Hence
\[
{\rm POL}\cap\Delta(S)\neq\varnothing.
\]

Since $S$ was arbitrary among all $r$-element subsets of $\mathbb N$,
it follows that ${\rm POL}$ is a $\Delta_r^*$-set.
\end{proof}

We need the following lemma in our next proof.
\begin{lem}\label{lem:IPstar-syndetic}
If $A$ is $IP_{r}^{*}$ then $A$ is multiplicatively syndetic.
More precisely, for any $y_{1},\dots,y_{r}\in\N$, put
$E=\FS(y_{i})_{i=1}^{r}$. Then
\[
\N=\bigcup_{t\in E} t^{-1}A. 
\]
\end{lem}

\begin{proof}[Proof of Lemma]
Fix $n\in\N$. Then
\[
nE=\left\lbrace nt:t\in E\right\rbrace=\left\lbrace n\sum_{i\in F}y_{i}:F\neq\emptyset\right\rbrace
     =\left\lbrace\sum_{i\in F} n y_{i}:F\neq\emptyset\right\rbrace
     =\FS(n y_{i})_{i=1}^{r}.
\]
Thus $nE$ is an $IP_{r}$-set (of length $r$).
Since $A$ is $IP_{r}^{*}$, $nE\cap A\neq\emptyset$.
Hence there is $t\in E$ with $nt\in A$, i.e. $n\in t^{-1}A$.
As $n$ was arbitrary, $\N=\bigcup_{t\in E}t^{-1}A$.
\end{proof}

We also need the following lemma to prove Theorem  \ref{B}.
\begin{lem}
\label{thm:POL-IPstar}
Let $k=k_2\leq50$ and let $p_1<p_2<\cdots<p_{k-1}$
be the first $k-1$ primes. Define $$M_k=1+p_1+p_1p_2+\cdots+p_1p_2\cdots p_{k-1}.$$
Then, for every choice of positive integers $x_1,\ldots,x_{M_k},$
we have $FS(x_1,\ldots,x_{M_k})\cap{\rm POL}\neq\varnothing.$
Consequently, ${\rm POL}$ is an $IP_{M_k}^*$-set.
\end{lem}

\begin{proof}
Define integers $L_k=1$ and, recursively for $1\leq i<k$, $L_i=p_iL_{i+1}+1.$
This recurrence formulae gives $$L_1
=
1+p_1+p_1p_2+\cdots+p_1p_2\cdots p_{k-1}
=
M_k.$$

Fix arbitrary positive integers $x_1,\ldots,x_{M_k}.$
Define the partial sums
\[
s_j=x_1+\cdots+x_j,
\qquad
1\leq j\leq M_k,
\]
and put $C_1=\{s_1,\ldots,s_{M_k}\}.$
Since all $x_j$ are positive, the numbers
$s_1,\ldots,s_{M_k}$ are distinct, and hence $|C_1|=M_k=L_1.$

We now construct distinct elements $b_1,\ldots,b_k$
and nonempty sets $C_1\supseteq C_2\supseteq\cdots\supseteq C_k$ inductively.

Suppose that $C_i$ has been constructed and satisfies $|C_i|\geq L_i.$
Among the $p_i$ residue classes modulo $p_i$, one contains at least $\frac{|C_i|}{p_i}$
elements of $C_i$. Choose such a residue class, say $h_i\pmod{p_i},$
and choose $b_i\in C_i$
with $b_i\equiv h_i\pmod{p_i}.$
Define
\[
C_{i+1}
=
\bigl(C_i\cap(h_i+p_i\mathbb Z)\bigr)\setminus\{b_i\}.
\]
Then
\[
|C_{i+1}|
\geq
\frac{|C_i|}{p_i}-1
\geq
\frac{L_i-1}{p_i}.
\]
Since $L_i=p_iL_{i+1}+1,$
we obtain $|C_{i+1}|\geq L_{i+1}.$
This completes the induction.

In particular, $C_k$ is nonempty, so we may choose $b_k\in C_k.$
Thus $b_1,\ldots,b_k$ are distinct elements of $C_1$.
We claim that $H=\{b_1,\ldots,b_k\}$ is admissible.

Fix $i\in\{1,\ldots,k-1\}$. By construction, $b_i\equiv h_i\pmod{p_i},$
and, since $C_{i+1}\subseteq h_i+p_i\mathbb Z,$
all later elements $b_{i+1},\ldots,b_k$
are also congruent to $h_i$ modulo $p_i$. Consequently, modulo $p_i$,
the set $H$ can occupy at most $i$ distinct residue classes.
The residue classes of $b_1,\ldots,b_{i-1}$
together with the single residue class $h_i$ occupied by $b_i,\ldots,b_k.$
For $i=1$ this gives exactly one residue class. For $i\geq2$, $p_i\geq i+1,$
and hence $|H\bmod p_i|\leq i<p_i.$

Now let $p\geq p_k$. Since $p_k>k$ and $|H|=k$, $|H\bmod p|\leq k<p.$
Therefore $H$ is admissible.

We next show that every difference of two elements of $H$ belongs to
the finite sumset $FS(x_1,\ldots,x_{M_k}).$
Indeed, every $b_i$ belongs to $C_1$, so there exist indices $1\leq u_i\leq M_k$
such that $b_i=s_{u_i}.$
If $b_i>b_j$, then necessarily $u_i>u_j$, and therefore
\[
b_i-b_j
=
s_{u_i}-s_{u_j}
=
x_{u_j+1}+\cdots+x_{u_i}
\in FS(x_1,\ldots,x_{M_k}).
\]
Thus $\Delta(H)
\subseteq
FS(x_1,\ldots,x_{M_k}).$

Since $|H|=k=k_2$ and $H$ is admissible, Theorem~\ref{ess}
with $m=2$ yields distinct $b_i,b_j\in H$ such that $n+b_i, n+b_j$
are consecutive primes for infinitely many $n\in\mathbb N$.
Hence $|b_i-b_j|\in{\rm POL}.$ 
Because $b_i,b_j$ belong to $C_1$ and the larger partial sum minus
the smaller partial sum is a finite sum of consecutive $x_\ell$'s,
we have
\[
|b_i-b_j|
\in
FS(x_1,\ldots,x_{M_k}).
\]
Therefore $FS(x_1,\ldots,x_{M_k})\cap{\rm POL}\neq\varnothing.$

Since the choice of $x_1,\ldots,x_{M_k}$ was arbitrary, ${\rm POL}$ is
an $IP_{M_k}^*$-set.
\end{proof}

\begin{proof}[Proof of Theorem \ref{B}]
From Lemma~\ref{thm:POL-IPstar}, we have ${\rm POL}$ is an
$IP_{M_k}^*$-set. In the Lemma~\ref{lem:IPstar-syndetic}, choose
\[
y_i=2^{i-1},
\qquad
1\leq i\leq M_k.
\]
Then every positive integer up to $2^{M_k}-1$ has a unique binary
expansion, and hence
\[
FS(1,2,4,\ldots,2^{M_k-1})
=
\{1,2,\ldots,2^{M_k}-1\}.
\]
Therefore, putting $E=\{1,2,\ldots,2^{M_k}-1\},$ from the 
Lemma~\ref{lem:IPstar-syndetic}, we have $\mathbb N = \bigcup_{t\in E}t^{-1}{\rm POL}.$ Moreover,  $|E|=2^{M_k}-1.$  Thus ${\rm POL}$ is multiplicatively syndetic.

Since $M_k = 1+p_1+p_1p_2+\cdots+p_1p_2\cdots p_{k-1},$ we have the elementary bound $M_k\leq k\prod_{i=1}^{k-1}p_i.$  Consequently,
\[
|E|
=
2^{M_k}-1
\leq
2^{\,k\prod_{i=1}^{k-1}p_i}.
\]
Since $k=k_2\leq50$, this gives the completely explicit bound $|E| \leq 2^{\,50\prod_{i=1}^{49}p_i}.$
This completes the proof.
\end{proof}

\subsection*{Acknowledgement} The author of this paper is supported by NBHM postdoctoral fellowship with reference no: 0204/27/(27)/2023/R \& D-II/11927.

\vspace{.5in}

\noindent
\textbf{Sayan Goswami}\\
Ramakrishna Mission Vivekananda Educational and Research Institute, Belur Math,\\
Howrah, West Bengal--711202, India.\\
Email: \href{mailto:sayan92m@gmail.com}{sayan92m@gmail.com}\\
Homepage: \url{https://sites.google.com/view/sayan-goswami}\\[1.5em]

\end{document}